\documentclass{article}
\usepackage{amsfonts}
\usepackage{amsmath}

\newtheorem{theorem}{Theorem}

\newtheorem{lemma}[theorem]{Lemma}

\newenvironment{proof}[1][Proof]{\noindent\textbf{#1.} }{\ \rule{0.5em}{0.5em}}
\input{tcilatex}

\begin{document}

\title{ Feedback Equivalence of $1$-dimensional Control Systems of the $1$%
-st Order}
\author{Valentin Lychagin \\
Department of Mathematics, University of Tromso\\
\&\\
Institute of Control Science , Russian Academy of Science}
\maketitle

\begin{abstract}
The problem of local feedback equivalence for $1$-dimensional control
systems of the $1$-st order is considered. The algebra of differential
invariants and criteria for the feedback equivalence for regular control
systems are found.
\end{abstract}

\section{Introduction}

In this paper we study the problem of local feedback equivalence for $1$%
-dimensional control systems of $1$-st order. 

As in paper (\cite{L}) we use the method of differential invariants. To this
end we consider control systems as underdetermined ordinary differential
equations. This gives a representation of feedback transformations as a
special type of Lie transformations, and we study and find differential
invariants of these representation.

Remark also that from the\ EDS point of view the case of control systems
considered here is equivalent to the case of second order systems considered
in (\cite{L}), but from ODE point of view they have different algebras of
feedback differential invariants.

To find a structure of the algebra of feedback differential invariants we
first  find $3$ feedback invariant derivations. Then the differential
invariants algebra is generated by two basic differential invariants $J$ and 
$K$ of orders $2$ and $3$ respectively and by all their invariant
derivations.

This description allows us to find invariants for the formal feedback
equivalence problem.

To get a local feedback equivalence we introduce a notion of \textit{regular 
}control system and connect with such a system a $3$-dimensional submanifold 
$\Sigma $ in $\mathbb{R}^{14}$\textit{. }

The main result of the paper states that two regular control systems are
locally feedback equivalent if and only if the corresponding $3$-dimensional
submanifolds $\Sigma $ coincide.

\section{ Representation of Feedback Pseudogroup}

Let 
\begin{equation}
\overset{\cdot }{x}=F(x,u,\overset{\cdot }{u}),  \label{Cont1}
\end{equation}%
be an autonomous 1-dimensional control system of the $1$-st order.

Here the function $x=x\left( t\right) $ describes a dynamic of the state of
the system, and $u=u\left( t\right) $ is a scalar control parameter.

We shall consider this system as an undetermined ordinary differential
equation of the first order on sections of $2$-dimensional bundle $\pi :%
\mathbb{R}^{3}\rightarrow \mathbb{R}$ , where $\pi :(x,u,t)\longmapsto t.$

Let $\mathcal{E}\subset J^{1}\left( \pi \right) $ be the corresponding
submanifold. In the canonical jet coordinates $\left(
t,x,u,x_{1},u_{1},....\right) $ this submanifold is given by the equation:%
\begin{equation*}
x_{1}=F\left( x,u,u_{1}\right) .
\end{equation*}

It is known (see, for example, \cite{KLV}) that Lie transformations in jet
bundles $J^{k}\left( \pi\right) $ for $2$-dimensional bundle $\pi$ are
prolongations of point transformations, that is, prolongations of
diffeomorphisms of the total space of the bundle $\pi.$

We shall restrict ourselves by point transformations which are automorphisms
of the bundle $\pi.$

Moreover, if these transformations preserve the class of systems (\ref{Cont1}%
) then they should have the form 
\begin{equation}
\Phi :\left( x,u,t\right) \rightarrow \left( X\left( x\right) ,U\left(
x,u\right) ,t\right) .  \label{Feed1}
\end{equation}

Diffeomorphisms of form (\ref{Feed1}) is called \emph{feedback
transformations. } The corresponding infinitesimal version of this notion is
a \emph{feedback vector field, }i.e. a plane vector field of the form%
\begin{equation*}
X_{a,b}=a\left( x\right) \partial _{x}+b\left( x,u\right) \partial _{u}.
\end{equation*}

The feedback transformations in a natural way act on the control systems of
type (\ref{Cont1}): 
\begin{equation*}
\mathcal{E}\longmapsto \Phi ^{\left( 1\right) }\left( \mathcal{E}\right) ,
\end{equation*}%
where $\Phi ^{\left( 1\right) }:J^{1}\left( \pi \right) \rightarrow
J^{1}\left( \pi \right) $ is the first prolongation of the point
transformation $\Phi .$

Passing to functions $F,$ defining the systems, we get the following action
on these functions:$\widehat{\Phi }:F\longmapsto G,$ where the function $G$
is a solution of the equation%
\begin{equation}
X_{x}~G=F\left( X,U,U_{x}G+U_{u}u_{1}\right) .  \label{rep1}
\end{equation}

The infinitesimal version of this action leads us to the following
representation $X_{a,b}\longmapsto \widehat{X_{a,b}}~$ of feedback vector
fields:

\begin{equation}
\widehat{X_{a,b}}~=a~\partial _{x}+b~\partial _{u}+\left(
u_{1}b_{u}+f~b_{x}\right) \partial _{u_{1}}+a_{x}~f~\partial _{f}.
\label{rep2}
\end{equation}

In this formula $\widehat{X_{a,b}}$ is a vector field on the $4$-dimensional
space $\mathbb{R}^{4}$ with coordinates $\left( u,u,u_{1},f\right) ,$ and
this field corresponds to the above action in the following sense. 

Each control system (\ref{Cont1}) determines a $3$-dimensional submanifold $%
L_{F}\subset \mathbb{R}^{4},$ the graph of $F:$%
\begin{equation*}
L_{F}=\left\{ f=F\left( x,u,u_{1}\right) \right\} .
\end{equation*}%
Let $A_{t}$ be the $1$-parameter group of shifts along vector field $X_{a,b}$
and let $B_{t}:\mathbb{R}^{4}\rightarrow \mathbb{R}^{4}$ be the
corresponding $1$-parameter group of shifts along $\widehat{X_{a,b}},$ then
these two actions related as follows 
\begin{equation*}
L_{\widehat{A_{t}}\left( F\right) }=B_{t}\left( L_{F}\right) .
\end{equation*}%
In other words, if we consider an $1$-dimensional bundle 
\begin{equation*}
\kappa :\mathbb{R}^{4}\rightarrow \mathbb{R}^{3},
\end{equation*}%
where $\kappa ((u,u,u_{1},f))=(u,u,u_{1}),$ then formula (\ref{rep2})
defines the representation $X\longmapsto \widehat{X}$ of the Lie algebra of
feedback vector fields into the Lie algebra of Lie vector fields on $%
J^{0}\left( \kappa \right) ,$ and the action of Lie vector fields $\widehat{X%
}$ on sections of bundle $\kappa $ corresponds to the action of feedback
vector fields on right hand sides of (\ref{Cont1})

\section{Feedback Differential Invariants}

By a \emph{feedback differential invariant} of order $\leq k$ we understand
a function $I\in C^{\infty}\left( J^{k}\kappa\right) $ on the space of $k$%
-jets $J^{k}(\kappa),$ which is invariant under of the prolonged action of
feedback transformations.

Namely, 
\begin{equation*}
\widehat{X_{a,b}}^{\left( k\right) }\left( I\right) =0,
\end{equation*}
for all feedback vector fields $X_{a,b}.$

In what follows we shall omit subscript of order of jet spaces, and say that
a function $I$ on the space of infinite jets $I\in C^{\infty }\left(
J^{\infty }\kappa \right) $ is a feedback differential invariant if 
\begin{equation*}
\widehat{X_{a,b}}^{\left( \cdot \right) }\left( I\right) =0,
\end{equation*}%
where $\widehat{X_{a,b}}^{\left( \cdot \right) }$ is the prolongation of the
vector field $X_{a,b}$ in the space of infinite jets $J^{\infty }\kappa .$

In a similar way one defines a \emph{feedback invariant derivations} as
combinations of total derivatives%
\begin{align*}
& \nabla =A\frac{d}{dx}+B\frac{d}{du}+C\frac{d}{du_{1}}, \\
& A,B,C\in C^{\infty }\left( J^{\infty }\kappa \right) ,
\end{align*}%
which are invariant with respect to prolongations of feedback
transformations, that is,%
\begin{equation*}
\lbrack \widehat{X_{a,b}}^{\left( \cdot \right) },\nabla ]=0
\end{equation*}%
for all feedback vector fields $X_{a,b}.$

Remark that for these derivations functions $\nabla \left( I\right) $ are
differential invariants ( of order, as a rule, higher then order of $I$) for
any feedback differential invariant $I.$ This observation allows us to
construct new differential invariants from known ones only by the
differentiations.

Recall the construction of the Tresse derivations in our case. Let $%
J_{1},J_{2},J_{3}\in C^{\infty }\left( J^{k}\kappa \right) $ be three
feedback differential invariants, and let 
\begin{equation*}
\widehat{d}J_{i}=\frac{dJ_{i}}{dx}dx+\frac{dJ_{i}}{du}du+\frac{dJ_{i}}{du_{1}%
}du_{1}
\end{equation*}

be their total derivatives.

Assume that we are in a domain $\mathcal{D}$ in $J^{k}\kappa ,$ where%
\begin{equation*}
\widehat{d}J_{1}\wedge \widehat{d}J_{2}\wedge \widehat{d}J_{3}\neq 0.
\end{equation*}%
Then, for any function $V\in C^{\infty }\left( J^{l}\kappa \right) $ over
domain $\mathcal{D}$, one has decomposition 
\begin{equation*}
\widehat{d}V=\lambda _{1}\widehat{d}J_{1}+\lambda _{2}\widehat{d}%
J_{2}+\lambda _{3}\widehat{d}J_{3}.
\end{equation*}%
Coefficients $\lambda _{1},\lambda _{2}$ and $\lambda _{3}$ of this
decomposition are called the \textit{Tresse derivatives} of $V$ and are
denoted by 
\begin{equation*}
\lambda _{i}=\frac{DV}{DJ_{i}}.
\end{equation*}%
The remarkable property of these derivatives is the fact that they are
feedback differential invariants (of higher, as a rule, order then $V$ )
each time when $V$ is a feedback differential invariant.

In other words, the Tresse derivatives 
\begin{equation*}
\frac{D}{DJ_{1}},\frac{D}{DJ_{2}}\ \text{and }\frac{D}{DJ_{3}}
\end{equation*}%
are feedback invariant derivations.

\section{Dimensions of Orbits}

First of all, we remark that the submanifold $\left\{ f=0\right\} $ is a
singular orbit for the feedback action in the space of $0$-jets $J^{0}\kappa 
$. The generating function of the feedback vector field $\widehat{X_{a,b}}~\ 
$has the form:%
\begin{equation*}
\phi _{a,b}=a_{x}f-af_{x}-bf_{u}-\left( u_{1}b_{u}+fb_{x}\right) f_{z},
\end{equation*}%
and the formula\ for prolongations of vector fields (\cite{KLV}) shows that
in the space of $1$-jets $J^{1}\kappa ,$ in addition, we have one more
singular orbit $\left\{ f_{u_{1}}=0\right\} .$ In similar way, we have one
more singular orbit $\left\{ f_{u_{1}u_{1}}=0\right\} $ in the space of $2$%
-jets. There are no more additional singular orbits in the spaces of $k$%
-jets, when $k\geq 3.$

We say that a point $x_{k}\in J^{k}\kappa $ is \emph{regular, } if $f\neq
0,f_{u_{1}}\neq 0,f_{u_{1}u_{1}}\neq 0$ at this point.

In what follows we shall consider orbits of regular points only.

It is easy to see, that the $k-$th prolongation of the feedback vector field 
$\widehat{X_{a,b}}$ depends on $\left( k+1\right) $-jet of function $a\left(
x\right) $ and ($k+1)$-jet of function $b\left( x,u\right) .$

Denote by $V_{i}^{k}$ and $W_{ij}^{k}$ the components of the decomposition%
\begin{equation*}
\widehat{X_{a,b}}^{\left( k\right) }=\dsum\limits_{0\leq i\leq k+1}a^{\left(
i\right) }\left( x\right) V_{i}^{k}+\dsum\limits_{0\leq i+j\leq k+1}\frac{%
\partial ^{i+j}b}{\partial x^{i}\partial u^{j}}W_{ij}^{k}.
\end{equation*}%
Then, by the construction, the vector fields $V_{i}^{k},0\leq i\leq k+1$,
and $W_{ij}^{k},0\leq i+j\leq k+1$, generate a completely integrable
distribution on the space of $k$-jets, integral manifolds of which are
orbits of the feedback action in $J^{k}\kappa .$

Straightforward computations show that there are no non trivial feedback
differential invariants of the $1$-st order.

Let $\mathcal{O}_{k+1}$ be a feedback orbit in $J^{k+1}\kappa $, then the
projection $\mathcal{O}_{k}=\kappa _{k+1,k}\left( \mathcal{O}_{k+1}\right)
\subset J^{k}\kappa $ is an orbit too, and to determine dimensions of the
orbits one should find dimensions of the bundles: $\kappa _{k+1,k}:\mathcal{O%
}_{k+1}\rightarrow \mathcal{O}_{k}.$ To do this we should find conditions on
functions $a$ and $b$ under which $\widehat{X_{a,b}}^{\left( k\right) }=0$
at a point $x_{k}\in J^{k}\kappa .$

Assume that $\widehat{X_{a,b}}^{\left( k-1\right) }=0$ at the point $%
x_{k-1}\in J^{k-1}\kappa $ . Then the vector field $\widehat{X_{a,b}}%
^{\left( k\right) }$ is a $\kappa _{k,k-1}$-vertical over this point. 

Components 
\begin{equation*}
\frac{d^{k}\phi }{dx^{i}du^{j}}\frac{\partial }{\partial f_{\sigma _{ij}}}
\end{equation*}%
of this vector field, where $\sigma _{ij}=(\underset{i\text{-times}}{%
\underbrace{x,....,x},}\underset{j\text{-times}}{\underbrace{u...,u}}),i+j=k,
$ and components 
\begin{equation*}
\frac{d^{k}\phi }{dx^{i}du^{j}du_{1}}\frac{\partial }{\partial f_{\tau _{ij}}%
},
\end{equation*}%
where $\tau _{ij}=(\underset{i\text{-times}}{\underbrace{x,....,x},}\underset%
{j\text{-times}}{\underbrace{u...,u}}),i+j=k-1$ depend on 
\begin{equation*}
\frac{\partial ^{k+1}b}{\partial x^{i}\partial u^{j}},
\end{equation*}%
and 
\begin{equation*}
\frac{d^{k+1}a}{dx^{k+1}}
\end{equation*}%
respectively.

All others components 
\begin{equation*}
\frac{d^{k}\phi }{dx^{r}du^{s}du_{1}^{t}}\frac{\partial }{\partial f_{\sigma
}}
\end{equation*}%
are expressed in terms of $k$-jet of $b\left( x,u\right) $ and $k$-jet of
function $a\left( x\right) .$

It shows that the bundles: $\kappa _{k,k-1}:\mathcal{O}_{k}\rightarrow 
\mathcal{O}_{k-1}$ are $(k+3)$-dimensional, when $k>1$.

Feedback orbits in the space of $2$-jets can be found by direct integration
of $12$-dimensional completely integrable distribution generating by the
vector fields $V_{i}^{1},0\leq i\leq 3$, and $W_{ij}^{1},0\leq i+j\leq 2$.
Summarizing, we get the following result.

\begin{theorem}
\begin{enumerate}
\item The first non-trivial differential invariants of feedback
transformations appear in order $2$ and they are functions of the basic
invariant 
\begin{equation*}
J=\frac{f^{2}~f_{u_{1}u_{1}}}{\left( u_{1}f_{u_{1}}-f\right) ~f_{u_{1}}^{2}}.
\end{equation*}

\item There are 
\begin{equation*}
\frac{k\left( k+1\right) }{2}-2
\end{equation*}%
independent differential invariants of pure order $k$.

\item Dimension of the algebra of differential feedback invariants of order $%
k\geq 2,$ is equal to%
\begin{equation*}
\frac{k^{3}}{6}+\frac{k^{2}}{2}-\frac{5k}{3}+1.
\end{equation*}

\item Dimension of the regular feedback orbits in the space \ of $k$-jets, $%
k\geq 2,$ is equal to 
\begin{equation*}
\frac{\left( k+1\right) ^{2}}{2}+\frac{23k}{3}+\frac{5}{2}.
\end{equation*}
\end{enumerate}
\end{theorem}

\section{Invariant Derivations}

We'll need the following result which allows us to compute invariant
derivations. 

Assume that an infinitesimal Lie pseudogroup $\mathfrak{g}$  is represented
in the Lie algebra of contact vector fields on the manifold of $1$-jets  $%
J^{1}\left( \mathbb{R}^{n}\right) .$ 

Moreover, we will identify elements $\mathfrak{g}$ with the corresponding
contact vector fields , i.e. we assume that elements of $\mathfrak{g}$ have
the form $X_{f}$ (see \cite{KLV}), where $f$ is the generating function.

\begin{lemma}
Let $x_{1},..,x_{n}$ be coordinates in $\mathbb{R}^{n}$, and let $\left(
x_{1},...,x_{n},u,p_{1},..,p_{n}\right) $ be the corresponding canonical
coordinates in the $1$-jet space $J^{1}\left( \mathbb{R}^{n}\right) .$ 

Then a derivation 
\begin{equation*}
\nabla =\sum_{i=1}^{n}A_{i}\frac{d}{dx_{i}}
\end{equation*}%
is $\mathfrak{g}$-invariant if and only if functions $A_{i}\in C^{\infty
}\left( J^{\infty }\mathbb{R}^{n}\right) ,\ j=1,..,n,$ are solutions of the
following PDE system:%
\begin{equation}
X_{f}\left( A_{i}\right) +\sum_{j=1}^{n}\frac{d}{dx_{j}}\left( \frac{%
\partial f}{\partial p_{i}}\right) A_{j}=0,  \label{InvariantDerSystemEq}
\end{equation}%
for all $i=1,...,n,$ and $X_{f}\in \mathfrak{g.}$
\end{lemma}

\begin{proof}
We have  (\cite{KLV}): 
\begin{equation*}
X_{f}^{\bullet }=\mathbf{E}_{f}-\sum_{i=1}^{n}\frac{\partial f}{\partial
p_{i}}\frac{d}{dx_{i}},
\end{equation*}%
where%
\begin{equation*}
\mathbf{E}_{f}=\sum_{\sigma }\frac{d^{\left\vert \sigma \right\vert }f}{%
dx^{\sigma }}\frac{\partial }{\partial p_{\sigma }}
\end{equation*}%
is the evolutionary derivation, $\sigma $ is a multi index and $\left\{
p_{\sigma }\right\} $ are the canonical coordinates in $J^{\infty }\mathbb{R}%
^{n}.$

Using the fact that evolutionary derivations commute with the total ones and
the relation 
\begin{equation*}
\lbrack \nabla ,X_{f}^{\bullet }]=0,
\end{equation*}%
we get%
\begin{eqnarray*}
0 &=&\left[ \sum_{j=1}^{n}A_{j}\frac{d}{dx_{j}},\mathbf{E}_{f}-\sum_{i=1}^{n}%
\frac{\partial f}{\partial p_{i}}\frac{d}{dx_{i}}\right]  \\
&=&-\sum_{j}\mathbf{E}_{f}\left( A_{j}\right) \frac{d}{dx_{j}}%
+\sum_{i,j}\left( -A_{j}\frac{d}{dx_{j}}\left( \frac{\partial f}{\partial
p_{i}}\right) \frac{d}{dx_{i}}+\frac{\partial f}{\partial p_{i}}\frac{dA_{j}%
}{dx_{i}}\frac{d}{dx_{j}}\right)  \\
&=&-\sum_{s}\left( X_{f}^{\bullet }\left( A_{s}\right) +\sum_{j}A_{j}\frac{d%
}{dx_{j}}\left( \frac{\partial f}{\partial p_{s}}\right) \right) \frac{d}{%
dx_{s}}.
\end{eqnarray*}
\end{proof}

In our case we expect three linear independent feedback invariant
derivations. To solve PDE system (\ref{InvariantDerSystemEq}) we first
assume that the unknown functions are functions on the $1$-jet space $J^{1}%
\mathbb{R}^{3}.$ Then collect terms in (\ref{InvariantDerSystemEq}) with $%
a,a^{^{\prime }},a^{\prime \prime }$ and $b,b_{x},b_{u},b_{xx},b_{xu}$ and $%
b_{u~u}$ we get the system of $8$ differential equations for $3$ unknown
functions. Solving the system we found two independent invariant
derivations. The last one we get in a similar way by assuming that the
unknown functions are functions on the $2$-jet space $J^{2}\mathbb{R}^{3}.$

Finally, we have $3$ feedback invariant derivations: 
\begin{eqnarray*}
\nabla _{1} &=&\frac{u_{1}f_{u_{1}}-f}{f_{u_{1}}}\frac{d}{du}+\frac{%
f-u_{1}f_{u_{1}}}{f_{u_{1}}^{2}}f_{u}\frac{d}{du_{1}}, \\
\nabla _{2} &=&\frac{f}{f_{u_{1}}}\frac{d}{du_{1}}, \\
\nabla _{3} &=&f~\frac{d}{dx}+\frac{f}{f_{u_{1}}}\frac{d}{du}+ \\
&&\left( \frac{f_{x}f_{u_{1}}+f_{u}-zf_{u~u_{1}}-f_{xu_{1}}}{f_{u_{1}u_{1}}}+%
\frac{u_{1}f_{u_{1}}-f}{f_{u_{1}}^{2}}f_{u}\right) \frac{d}{du_{1}}.
\end{eqnarray*}%
These derivations obey the following commutation relations%
\begin{eqnarray*}
\lbrack \nabla _{2},\nabla _{1}] &=&J~\nabla _{1} \\
\lbrack \nabla _{3},\nabla _{1}] &=&K~\nabla _{2} \\
\lbrack \nabla _{3},\nabla _{2}] &=&-\nabla _{3}+J~\nabla _{1}+L~\nabla _{2}
\end{eqnarray*}%
where $K$ and $L$ are some differential invariants of the 3rd order (see
below).

\bigskip

\section{Differential Invariants of the $3$-rd Order}

Theorem 1 shows that there are four independent differential invariants of
the $3$-rd order. We get three of them simply by invariant differentiations: 
\begin{equation*}
\nabla _{1}\left( J\right) ,\nabla _{2}\left( J\right) ,\nabla _{3}\left(
J\right) .
\end{equation*}%
The symbols of these invariants contain:

\begin{itemize}
\item symbol of $\nabla _{2}\left( J\right) $ depends on $%
f_{u_{1}u_{1}u_{1}},$

\item symbol of $\nabla _{1}\left( J\right) $ depends on $%
f_{u_{1}u_{1}u_{1}} $ and $f_{uu_{1}u_{1}},$

\item symbol of $\nabla _{3}\left( J\right) $ depends on $%
f_{u_{1}u_{1}u_{1}} $,$f_{uu_{1}u_{1}}$ and $f_{xu_{1}u_{1}}.$
\end{itemize}

It shows that these differential invariants are independent.

The similar observation shows that the differential invariant $L,$ which
appears in the commutation relations,  is a function of $J,\nabla _{1}\left(
J\right) ,\nabla _{2}\left( J\right) ,\nabla _{3}\left( J\right) ,$ and the
differential invariant $K$ is the forth independent invariant.  It has the
following form:%
\begin{eqnarray*}
K &=&-u_{1}f_{xu}+2u_{1}{\frac{\,f_{u}^{2}}{ff_{u_{1}}}}-2\frac{\,f_{u}^{2}}{%
{f_{u_{1}}}^{2}} \\
&&+{\frac{f_{~uu}u_{1}-2\,f_{u}f_{x}+ff_{xu}}{f_{u_{1}}}}-u_{1}{\frac{\left(
f_{~uu}u_{1}-2\,f_{u}f_{x}\right) }{f}} \\
&&+{\frac{c_{1}}{f_{u_{1}}{f_{u_{1}u_{1}}}^{2}}+\frac{c_{2}}{f{f_{u_{1}u_{1}}%
}^{2}}+\frac{c_{3}}{{f_{u_{1}u_{1}}}^{2}}}+{\frac{c_{4}}{ff_{u_{1}u_{1}}}}+{%
\frac{c_{5}}{f_{u_{1}}f_{u_{1}u_{1}}}}+{\frac{c_{6}}{f_{u_{1}u_{1}}},}
\end{eqnarray*}%
where

\begin{eqnarray*}
&&c_{1}=-ff_{u}f_{xu_{1}}f_{u_{1}u_{1}u_{1}}-u_{1}f_{u}f_{uu_{1}}f_{u_{1}u_{1}u_{1}}+f_{u}^{2}f_{u_{1}u_{1}u_{1}},
\\
&&c_{2}=u_{1}\left(
f_{u}f_{u_{1}}f_{uu_{1}u_{1}}-f_{u}^{2}f_{u_{1}u_{1}u_{1}}-f_{x}f_{u}f_{u_{1}}f_{u_{1}u_{1}u_{1}}+f_{x}%
{f_{u_{1}}}^{2}f_{uu_{1}u_{1}}\right) \\
&&+{u_{1}}^{2}f_{uu_{1}}\left(
-f_{u_{1}}f_{uu_{1}u_{1}}+f_{u}f_{u_{1}u_{1}u_{1}}\right) , \\
&&c_{3}=ff_{xu_{1}}f_{uu_{1}u_{1}}+f_{x}f_{u}f_{u_{1}u_{1}u_{1}}-f_{u}f_{uu_{1}u_{1}}-f_{x}f_{u_{1}}f_{uu_{1}u_{1}}
\\
&&+u_{1}\left(
f_{u}f_{xu_{1}}f_{u_{1}u_{1}u_{1}}-f_{u_{1}}f_{xu_{1}}f_{uu_{1}u_{1}}+f_{uu_{1}}f_{uu_{1}u_{1}}\right) ,
\\
&&c_{4}=-u_{1}\left(
2\,f_{u_{1}}f_{x}f_{uu_{1}}-f_{u_{1}}f_{u}f_{xu_{1}}+f_{u}f_{uu_{1}}+f_{u_{1}}f_{~uu}+%
{f_{u_{1}}}^{2}f_{xu}\right) \\
&&+{u_{1}}^{2}\left( f_{u_{1}}f_{~uuu_{1}}-f_{u}f_{uu_{1}u_{1}}+{f_{uu_{1}}}%
^{2}\right) , \\
&&c_{5}=ff_{u}f_{xu_{1}u_{1}}-ff_{xu_{1}}f_{uu_{1}}+f_{u}f_{uu_{1}}+u_{1}%
\left( f_{u}f_{uu_{1}u_{1}}-{f_{uu_{1}}}^{2}\right) , \\
&&c_{6}=f_{~uu}-f_{u}f_{xu_{1}}+2%
\,f_{x}f_{uu_{1}}+f_{u_{1}}f_{xu}-ff_{~xuu_{1}} \\
&&+u_{1}\left(
f_{u_{1}}f_{~xuu_{1}}-f_{~uuu_{1}}+f_{xu_{1}}f_{uu_{1}}-f_{u}f_{xu_{1}u_{1}}%
\right) .
\end{eqnarray*}

\section{Algebra of Feedback Differential Invariants}

By \textit{regular orbits} we mean feedback orbits of regular points.

Counting the dimensions of regular feedback orbits shows that the following
result is valid.

\begin{theorem}
\label{MainDiffInv}Algebra of feedback differential invariants in a
neighborhood of a regular orbit is generated by differential invariant $J$
of the $2$-nd order, differential invariant $K$ of the $3$-rd order and all
their invariant derivatives.
\end{theorem}

\section{The Feedback Equivalence Problem}

Consider two control systems given by functions $F$ and $G.$ Then, to
establish feedback equivalence, we should solve the differential equation%
\begin{equation}
F\left( X,U,U_{x}G\left( x,u,u_{1}\right) +U_{u}u_{1}\right) -X_{x}~G\left(
x,u,u_{1}\right) =0  \label{equiveq1}
\end{equation}%
with respect to unknown functions $X\left( x\right) $ and $U\left(
x,u\right) .$

Let us denote the left hand side of (\ref{equiveq1}) by $H.$ Then assuming
the general position one can find functions $X,X_{x},U,U_{x},U_{u}$ from the
equations%
\begin{equation*}
H=H_{u_{1}}=H_{u_{1}}^{\left( 2\right) }=H_{u_{1}}^{\left( 3\right)
}=H_{u_{1}}^{\left( 4\right) }=0.
\end{equation*}

Remark, that the above general conditions are feedback invariant, depends on
finite jet of the system and holds in a dense open domain of the jet space.
Therefore, it holds in regular points.

Assume that we get%
\begin{align*}
U& =A\left( x,u,u_{1}\right) ,U_{x}=B\left( x,u,u_{1}\right) , \\
U_{u}& =C\left( x,u,u_{1}\right) ,X=D\left( x,u,u_{1}\right) , \\
X^{\prime }& =E\left( x,u,u_{1}\right)
\end{align*}%
Then the conditions 
\begin{align*}
A_{u_{1}}& =B_{u_{1}}=C_{u_{1}}=D_{u_{1}}=E_{u_{1}}=0, \\
D_{u}& =E_{u}=0
\end{align*}%
and 
\begin{equation*}
B=A_{x},C=A_{u},E=D_{x}
\end{equation*}%
show that if (\ref{equiveq1}) has a formal solution at each point $\left(
x,u,u_{1}\right) $ in a domain then this equation has a local smooth
solution.

On the other hand if system $F$ at a point $p=(x^{0},u^{0},u_{1}^{0})$ and
system $G$ at a point $\widetilde{p}=(\widetilde{x}^{0},\widetilde{u}^{0},%
\widetilde{u}_{1}^{0})$ has the same differential invariants then, by the
definition, there is a formal feedback transformation which send the
infinite jet of $F$ at the point $p$ to the infinite jet of $G$ at the point 
$\widetilde{p}.$

Keeping in mind these observations and results of theorem \ref{MainDiffInv}
we consider the space $\mathbb{R}^{3}$with coordinates $\left(
x,u,u_{1}\right) $ and the space $\mathbb{R}^{14}$ with coordinates $\left(
j,j_{1},j_{2},j_{3},j_{11},j_{12},j_{13},j_{22},j_{23},j_{33}k,k_{1},k_{2},k_{3}\right) . 
$

Then any control system, given by  function $F\left( x,u,u_{1}\right) $,
defines a map 
\begin{equation*}
\sigma _{F}:\mathbb{R}^{3}\rightarrow \mathbb{R}^{14},
\end{equation*}%
by 
\begin{eqnarray*}
j &=&J^{F},k=K^{F}, \\
j_{i} &=&\left( \nabla _{i}(J)\right) ^{F},k_{i}=\left( \nabla
_{i}(K)\right) ^{F}, \\
j_{ij} &=&\left( \nabla _{i}\nabla _{j}(J)\right) ^{F},
\end{eqnarray*}%
where $i,j=1,2,3,$ and the subscript $F$ means that the differential
invariants are evaluated due to the system.

Let 
\begin{equation*}
\Phi :\mathbb{R}^{3}\rightarrow \mathbb{R}^{3}
\end{equation*}%
be a feedback transformation.

Then from the definition of the feedback differential invariants it follows
that 
\begin{equation*}
\sigma _{F}\circ \Phi =\sigma _{\widehat{\Phi }\left( F\right) }.
\end{equation*}%
Therefore, the geometrical image 
\begin{equation*}
\Sigma _{F}={Im}\left( \sigma _{F}\right) \subset \mathbb{R}^{14}
\end{equation*}%
does depend on the feedback equivalence class of $F$ only.

We say that a system $F$ is \emph{regular} in a domain $D\subset\mathbb{R}%
^{3}$ if

\begin{enumerate}
\item $4$-jets of $F$ belong to regular orbits,

\item $\sigma _{F}\left( D\right) $ is a smooth $3$-dimensional submanifold
in $\mathbb{R}^{14},$ and

\item three of five functions $j,j_{1},j_{2},j_{3},k$ are coordinates on $%
\Sigma _{F}.$
\end{enumerate}

Assume, for example, that functions $j_{1},j_{2},j_{3}$ are coordinates on $%
\Sigma _{F}.$ The following lemma gives a relation between the Tresse
derivatives and invariant differentiations $\nabla _{1},\nabla _{2},\nabla
_{3}$.

\begin{lemma}
Let 
\begin{equation*}
\frac{D}{DJ_{1}},\frac{D}{DJ_{2}},\frac{D}{DJ_{3}}
\end{equation*}%
be the Tresse derivatives with respect to differential invariants $%
J_{i}=\nabla _{i}\left( J\right) .$

Then the following decomposition%
\begin{equation}
\nabla _{i}=\sum_{j}R_{ij}\frac{D}{DJ_{j}}  \label{invDecomp}
\end{equation}%
with feedback differential invariants $R_{ij}$ of order $\leq 4$ is valid.
\end{lemma}

\begin{proof}
Applying both parts of (\ref{invDecomp}) to invariant $J_{k}$ we get 
\begin{equation*}
\nabla _{i}\left( J_{k}\right) =R_{ik}
\end{equation*}%
which is a feedback differential invariant of order $\leq 4.$
\end{proof}

\begin{theorem}
Two regular systems $F$ and $G$ are locally feedback equivalent if and only
if 
\begin{equation}
\Sigma_{F}=\Sigma_{G}.  \label{EquivCond}
\end{equation}
\end{theorem}

\begin{proof}
Let us show that the condition \ref{EquivCond} implies a local feedback
equivalence.

Assume that 
\begin{align*}
J^{F}& =j^{F}\left( J_{1},J_{2},J_{3}\right) ,J_{ij}^{F}=j_{ij}^{F}\left(
J_{1},J_{2},J_{3}\right) , \\
K^{F}& =k^{F}\left( J_{1},J_{2},J_{3}\right) ,K_{i}^{F}=k_{i}^{F}\left(
J_{1},J_{2},J_{3}\right)
\end{align*}%
on $\Sigma _{F},$ and%
\begin{align*}
J^{G}& =j^{G}\left( J_{1},J_{2},J_{3}\right) ,J_{ij}^{G}=j_{ij}^{G}\left(
J_{1},J_{2},J_{3}\right) , \\
K^{G}& =k^{G}\left( J_{1},J_{2},J_{3}\right) ,K_{i}^{G}=k_{i}^{G}\left(
J_{1},J_{2},J_{3}\right)
\end{align*}%
on $\Sigma _{G}.$

Then condition \ref{EquivCond} shows that $%
j^{F}=j^{G},j_{ij}^{F}=j_{ij}^{G},k_{i}^{F}=k_{i}^{G}$ and $k^{F}=k^{G}.$

Moreover,as we have seen the invariant derivations $\nabla _{1},\nabla
_{2},\nabla _{3}$ are linear combinations of the Tresse derivatives with
coefficients which are feedback differential invariants of order $\leq 4$.

In other words, the above functions $j^{F},k^{F},j_{ij}^{F},k_{i}^{F}$ and
their partial derivatives in $j_{1},j_{2},j_{3}$ determine the restrictions
of all differential invariants.

Therefore, condition \ref{EquivCond} equalize restrictions of differential
invariants not only to order $\leq 4$ but in all orders, and provides 
formal and therefore local feedback equivalence between $F$ and $G$.
\end{proof}

\end{document}